\documentclass[11 pt]{amsart}
\usepackage{amsmath, amsthm, amssymb}
 \usepackage{amsfonts}
\usepackage{amsmath}
\usepackage{amsbsy}
\usepackage[all,cmtip]{xy}
\usepackage{amscd}
\usepackage{enumitem}
\usepackage[all]{xy}
\usepackage{graphicx}
\usepackage{ifpdf}
\usepackage[UKenglish]{babel}
\usepackage{hyperref}
\usepackage{color}

\newcommand{\Z}{\mathbb{Z}}
\newcommand{\F}{\mathbb{F}}

\newcommand{\tc}{\textcolor}

\newcommand{\da}{\ar@{-->}}
\newcommand{\dar}{\ar@{.>}}
\newcommand{\lar}{\ar@{-}}
\newcommand{\noar}{\ar@{}[]}
\newcommand{\del}{\partial}

\newcommand{\mf}{\mathfrak}

\newcommand{\bu}{\bullet}

\newtheorem{theorem}{Theorem}[section]

\newtheorem{prop}[theorem]{Proposition}
\newtheorem{corollary}[theorem]{Corollary}

\newtheorem{remark}[theorem]{Remark}

\begin{document}

\author{David Krcatovich}
\title[Alexander polynomials of $L$-space knots.]{A restriction on the Alexander polynomials of $L$-space knots.}
\date{\vspace{-3ex}}

\begin{abstract}
Using an invariant defined by Rasmussen, we extend an argument given by Hedden and Watson which further restricts which Alexander polynomials can be realized by $L$-space knots.
\end{abstract}

\maketitle

\section{Introduction}
In \cite{OSinteger}, Ozsv\'ath and Szab\'o show how the filtered chain homotopy type of the knot Floer complex $CFK^-(K)$ can be used to compute the Heegaard Floer homology of $S^3_n(K)$, the rational homology sphere obtained by doing Dehn surgery along $K\subset S^3$ with slope $n$. In \cite{OSlens}, they use this relationship to investigate which knots admit lens space surgeries, using the fact that if $Y$ is a lens space, $Y$ has the ``smallest possible" Heegaard Floer homology;
\begin{equation}
\text{rank}\ \widehat{HF}(Y) = | H_1(Y;\Z) |.
\label{simplehf}
\end{equation}
More generally, a rational homology sphere which satisfies condition \eqref{simplehf} is called an $L$-{\it{space}}. So, from a Heegaard-Floer perspective, a natural extension of the question ``which knots admit lens space surgeries?" is ``which knots admit $L$-space surgeries?".

Letting $A(x)$ denote the Alexander grading of a homogeneous element $x$ in $CFK^-$, the following Proposition is a straightforward consequence of \cite[Theorem 1.2]{OSlens} (cf. \cite[Remark 6.6]{HomConcordance}).
\begin{prop}
Suppose $K\subset S^3$ is a knot on which some positive integral surgery yields an $L$-space. Then $CFK^-(K)$ has a basis $\{ x_{-k}, \ldots, x_k \}$ with the following properties:
\begin{list}{$\bullet$}{}  
\item $A(x_i)=n_i$, where $n_{-k}<n_{-k+1} <\cdots < n_{k-1}<n_k$ 
\item $n_i = -n_{-i}$
\item If $i \equiv k \mod 2$, then $\del(x_i)=0$
\item If $i \equiv k+1 \mod 2$, then $\del(x_i) = x_{i-1} + U^{n_{i+1}-n_i}x_{i+1}$ \qed
\end{list}
\label{lspacechar}
\end{prop}
Notice that $x_k$ is in the kernel of $\del$ and not in the image, so
\begin{equation}
x_k \text{ generates }  H_* \left(\widehat{CFK}(K),\del\right) \cong \widehat{HF}(S^3) \cong \F .
\label{generator}
\end{equation}
By convention then, $M(x_k)=0$ (where $M$ is the Maslov grading). Since $U$ decreases $M$ by 2, and $\del$ decreases $M$ by 1, this determines the Maslov grading on all homogeneous elements of $CFK^-(K)$. 

Ozsv\'ath and Szab\'o also showed \cite{OSknot} that the graded Euler characteristic of $\widehat{CFK}(K)$ is the symmetrized Alexander polynomial of $K$,
\begin{equation}
\sum_i \chi \left(\widehat{CFK}(K,i)\right)\cdot T^i = \Delta_K(T),
\label{eulerhat}
\end{equation}
so a corollary to Proposition \ref{lspacechar} is
\begin{corollary}\cite[Corollary 1.3]{OSlens}
If $K\subset S^3$ is a knot which admits an $L$-space surgery, then
\begin{equation}
\Delta_K(T) = \sum_{i=-k}^{k} (-1)^{k+i}\ T^{n_i}
\label{OSpoly}
\end{equation}
for some sequence of integers \( n_{-k}<n_{-k+1} <\cdots < n_{k-1}<n_k \) satisfying \(n_{-i}=n_i\).
\label{OSrestriction}
\end{corollary}
\begin{remark}
Although Proposition \ref{lspacechar} only applies to knots which have positive $L$-space surgeries, a knot $K$ has a negative $L$-space surgery if and only if its mirror image $\overline{K}$ has a positive $L$-space surgery. Corollary \ref{OSrestriction} then follows in this generality because $\Delta_K(T)=\Delta_{\overline{K}}(T)$.
\label{mirror}
\end{remark}
In particular, all of the nonzero coefficients of $\Delta_K(T)$ are $\pm 1$. Note that 
\begin{equation}
n_k=g(K) = |\tau(K)|= g_4(K),
\label{genera}
\end{equation}
where $g$ is the Seifert genus, $\tau$ is the Ozsv\'ath-Szab\'o concordance invariant defined in \cite{OSfourballgenus}, and $g_4$ is the smooth four-genus of $K$. The first equality follows from the knot Floer homology detection of genus \cite{OSGenusBounds}, the second follows from \eqref{generator}, and the third follows from the fact shown in \cite{OSfourballgenus} that, for any knot $K$,\[ |\tau(K)|\leq g_4(K) \leq g(K).\]

This was the most general restriction on the Alexander polynomials of knots admitting $L$-space surgeries in the literature until Hedden and Watson showed\footnote{This result was given as a corollary of a more general restriction on knot Floer complexes, and was evidently already known to Rasmussen.}
\begin{prop}\cite[Corollary 9]{HeddenWatson}
If $K\subset S^3$ is a knot which admits an $L$-space surgery, then $\Delta_K(T)$ is as described in Corollary \ref{OSrestriction}, and further, $n_k - n_{k-1}=1$.
\label{HWrestriction}
\end{prop}
Their proof hinges on an invariant defined by Rasmussen, and a particular inequality satisifed by this invariant. Roughly, if large $n$-surgery is done on an unknot and on a knot $K$, the differences in the $d$-invariants (defined in Equation \eqref{ddef}) of the resulting manifolds are bounded above by numbers depending on $g_4(K)$. Proposition \ref{HWrestriction} is then proved by showing that if a complex has a basis as in Proposition \ref{lspacechar} and $n_k-n_{k-1} >1,$ Rasmussen's inequality is violated, and therefore this complex cannot be the knot Floer complex of any knot.

Our aim here is to extend this argument. We will introduce Rasmussen's invariant and inequality in Section \ref{invt}. In Section \ref{new}, we will show how to compute the invariant for $L$-space knots from their Alexander polynomials (that is, from the sequence of $n_i$'s). We will then see that Rasmussen's inequality places further restrictions on the $n_i$'s, analogous to the restriction $n_k-n_{k-1}=1$. As a result, it will be shown that certain symmetric Laurent polynomials satisfying Proposition \ref{HWrestriction} cannot be the Alexander polynomial of any $L$-space knot. 

\begin{theorem}
Suppose $K\subset S^3$ is a knot which admits an $L$-space surgery. Then its symmetrized Alexander polynomial can be written as \[\Delta_K(T) = \sum_{i=-k}^{k} (-1)^{k+i}\ T^{n_i},\] for some sequence of integers $n_{-k}<n_{-k+1} <\cdots < n_{k-1}<n_k$  satisfying the following:
\begin{list}{$\circ$}{}  
\item $n_i = -n_{-i}$
\item  if we let $r_i=n_{k+2-2i}-n_{k+1-2i}$, then $r_1=1,$ and for any $j\leq k$,   
\begin{equation}
 \sum_{i=2}^{j} r_i \  \leq \  \sum_{i=k-j+2}^{k} r_i.
\label{riineq}
\end{equation}
\end{list}
\label{newrestriction}
\end{theorem}
As we will explain in Section \ref{new}, the restriction is more easily stated in terms of a modified version of the Alexander polynomial, \[\widetilde{\Delta}_K(T):=\frac{\Delta_K(T)}{1-T^{-1}}.\] It follows from Corollary \ref{OSrestriction} that when $K$ is a knot which admits an $L$-space surgery, \[\widetilde{\Delta}_K(T)=\sum_{i=0}^\infty T^{a_i},\] for some sequence of integers satisfying 
\begin{list}{$\circ$}{}  
\item $a_0=g(K)$ 
\item $a_{i+1} < a_i$ 
\item $a_i=-i$ for $i\geq g(K).$
\end{list}
We can then rephrase Theorem \ref{newrestriction} as

\begin{theorem}[\sc Restatement of Theorem \ref{newrestriction} in terms of $\widetilde{\Delta}$ ]
Suppose $K\subset S^3$ is a knot which admits an $L$-space surgery and $\{a_i\}$ is the decreasing sequence of integers such that \[\widetilde{\Delta}_K(T) = \sum_{i=0}^{\infty}T^{a_i}.\] Then, for all $0\leq i \leq g(K)$, 
\begin{equation}
a_i\leq g(K)-2i.
\label{aiineq}
\end{equation}
\label{alternate}
\end{theorem}
To see the preceding Theorems as generalizations of Proposition \ref{HWrestriction}, note that in the language of Theorem \ref{newrestriction}, Proposition \ref{HWrestriction} translates to the statement $r_1=1$; in the language of Theorem \ref{alternate}, it translates to $a_1 \leq g(K)-2.$ \\

\noindent {\bf Acknowledgments.} The author would like to thank his advisor, Matt Hedden, for explaining the proof of Proposition \ref{HWrestriction}, and for suggesting the restatement of Theorem \ref{newrestriction} in terms of the polynomial $\widetilde{\Delta}$. The author was partially supported by NSF grant DMS--1150872.

\section{The invariant $h_m(K)$}
\label{invt}
A useful feature of Heegaard Floer theory is that its groups satisfy surgery exact triangles; for example, a long exact sequence between Heegaard Floer homology groups of manifolds which are $0-$, $\infty-$ and $n-$framed surgery along the same knot $K$ \cite[Section 9]{OSapplications}. In \cite[Definition 7.1]{Rasmussenknot}, Rasmussen defines an invariant $h_m(K)$ as the rank of a particular map in such a sequence (cf. \cite{Froyshovhdef}, where Fr{\o}yshov introduces an instanton-Floer invariant $h$).

Recall that if $(Y,\mf t)$ is a spin$^c$ rational homology sphere, Ozsv\'ath and Szab\'o define the $d$-{\it{invariant}} of $(Y,\mf t)$ as
\begin{equation}
d(Y,\mf t) = \min \{ M(x) | x\in \text{Im} \left( \pi_* :HF^\infty(Y,\mf t) \to HF^+(Y,\mf t)\right) \}.
\label{ddef}
\end{equation}
In \cite[Section 2.2]{Rasmussenh}, Rasmussen shows that, in the case where $S^3_{-n}(K)$ is an $L$-space, the invariant $d(S^3_{-n}(K), \mf{s}_m)$ is equal to twice $h_m(K)$, up to a shift which is independent of $K$. In particular, since $h_m({\text{unknot}})=0$ for all $m$, we have
\begin{equation}
h_m(K) = \frac{ d\left(S^3_{-n}(K),\mf{s}_m\right) - d\left(S^3_{-n}({\text{unknot}}),\mf{s}_m\right)}{2}.
\label{hdef}
\end{equation}
The key to obtaining restrictions on $L$-space knots is the following inequality, analagous to an inequality in instanton Floer homology proved by Fr{\o}yshov \cite{Froyshov}.

\begin{prop}\cite[Theorem 2.3]{Rasmussenh}
Let $K$ be a knot in $S^3$ and let $g_4(K)$ be its slice genus. Then $h_m(K)=0$ for $|m|>g_4(K)$, while for $|m|\leq g_4(K)$,
\begin{equation}
h_m(K)\leq \left\lceil \frac{g_4(K)-|m|}{2} \right\rceil.
\label{hineq}
\end{equation}
\label{hprop}
\end{prop}

Note that for a knot admitting an $L$-space surgery, due to \eqref{genera}, we can replace $g_4(K)$ with $g(K)$ and obtain
\begin{equation}
h_m(K)\leq \left\lceil \frac{g(K)-|m|}{2} \right\rceil.
\label{lspacehineq}
\end{equation}

It will be convenient to consider $L$-space knots -- that is, knots with {\it positive}, rather than negative, $L$-space surgeries. This is opposite Rasmussen's point of view in \cite{Rasmussenh}, but note that $K$ admits a positive $L$-space surgery if and only if its mirror image $\overline{K}$ admits a negative $L$-space surgery. Accordingly, we will follow Hedden and Watson in defining
\begin{equation}
\overline{h}_m(K):= \frac{d(S^3_n(\text{unknot}),\mf{s}_m) - d(S^3_n(K),\mf{s}_m)}{2}
\label{hbardef}
\end{equation}
and recall their observation that \(\overline{h}_m(K) = h_m(\overline{K}).\) Finally, we should note that $g(K)=g(\overline{K})$, so $\overline{h}_m$ satisfies the same inequality which $h_m$ does for knots admitting $L$-space surgeries; for $|m|\leq g(K)$,
\begin{equation}
\overline{h}_m(K)\leq \left\lceil \frac{g(K)-|m|}{2} \right\rceil.
\label{lspacehbarineq}
\end{equation}

\section{Values of $\overline{h}_m$ for $L$-space knots}
\label{new}

%
%
Next, we recall how to compute $d$-invariants, and therefore $\overline{h}_m$, from $CFK^-$. It was shown independently by Ozsv\'ath-Szab\'o \cite{OSknot} and Rasmussen \cite{Rasmussenknot} that for large $n$-surgery (that is, for $n\geq 2g(K)-1$), the Heegaard Floer homology groups $HF^-(S^3_n(K))$ are the homology groups of certain subcomplexes of $CFK^-(K)$, up to a shift in Maslov grading which is independent of $K$. In particular, if we let $A_m$ denote the subcomplex consisting of elements with Alexander grading less than or equal to $m$, then \[ HF^-(S^3_n(K), \mf{s}_m) \cong H_*(A_m),\] up to a shift in grading\footnote{Here we are adopting the convention that $CF^-$ and $CFK^-$ contain the element 1 in $\F[U]$.}. It follows that \[ d(S^3_n(K),\mf{s}_m) =\max \{ M(x) | x  \text{ a non-torsion generator of } H_*(A_m) \} +c,\] where $c$ is a constant which depends on $n$, but not on $K$. Therefore, the \textquotedblleft shifted" $d$-invariant  
\begin{equation}
 \widetilde{d}(K,m):=\max \{M(x) | x \text{ a non-torsion generator of } H_*(A_m) \},
\label{shiftddef}
\end{equation}
is well-defined, and satisfies
\begin{equation}
d(S^3_n(U),\mf s_m) - d(S^3_n(K),\mf s_m) = \widetilde{d}(U,m)-\widetilde{d}(K,m),
\label{shiftd}
\end{equation}
for any sufficiently large $n$. 
For the unknot, we have the complex \[CFK^- (\text{unknot}) \cong \F[U],\] where the generator has Maslov grading and Alexander grading equal to zero. Since multiplication by $U$ lowers the Alexander grading by 1 and the Maslov grading by 2, \[ \widetilde{d}(U,m)=m-|m|.\] Therefore, we can rewrite inequality \eqref{lspacehbarineq} using \eqref{shiftd} and the above: if $K \subset S^3$ is an $L$-space knot, then for $|m|\leq g(K)$,
\begin{equation}
\begin{split}
\overline{h}_m(K) = \frac{ \widetilde{d}(U,m)-\widetilde{d}(K,m)}{2} \leq &\ \left\lceil\frac{g(K)-|m|}{2}\right\rceil  \\
-\frac{1}{2}\widetilde{d}(K,m)\leq & \ \left\lceil\frac{g(K)-m}{2}\right\rceil.
\end{split}
\label{shiftdineq}
\end{equation}
%
%

With inequality \eqref{shiftdineq} in hand, it remains to see how the values of $\widetilde{d}$ are determined by the Alexander polynomial of an $L$-space knot.

%
%

Recall that the Alexander polynomial is the graded Euler characteristic of \( \widehat{CFK}(K)\), \[ \sum_i \chi \left( \widehat{CFK}(K,i)\right) \cdot T^i = \Delta_K(T).\] Further, $CFK^-$ is generated by the same set as $\widehat{CFK}$, over $\F[U]$ rather than $\F$. Since $U$ lowers the Alexander grading by 1 and preserves the parity of the Maslov grading, 
\begin{equation}
\begin{split}
& \ \ \ \ \sum_i \chi \left( CFK^-(K,i)\right) \cdot T^i \\
& = \  \sum_i \chi \left( \widehat{CFK}(K,i)\right) \cdot T^i \cdot (1+T^{-1} +T^{-2} +\cdots ) \\
& = \ \frac{\Delta_K(T)}{1-T^{-1}} \\
& =: \ \widetilde{\Delta}_K(T).
\end{split}
\end{equation}
In other words, \(\widetilde{\Delta}_K(T)\) is the graded Euler characteristic of $CFK^-(K)$. 
\begin{remark}
If $K\subset S^3$ is a knot for which \(\Delta_K(T)\) is of the type described in Corollary \ref{OSrestriction}, then \[\widetilde{\Delta}_K(T)= \sum_{i=0}^{\infty} T^{a_i},\] where
\begin{list}{$\bullet$}{}
\item \(a_0=g(K),\)
\item \( a_{i+1}<a_i,\) and 
\item \(a_i= -i\) for all \(i\geq g(K).\)
\end{list}
\label{lspacemodpoly}
\end{remark}

%
%

In \cite{Kreduced}, a reduced complex $\underline{CFK}^-$ was defined, and it was shown that for an $L$-space knot,
\begin{equation}
\underline{CFK}^-(K) \cong \F[U],
\label{lspacereduced}
\end{equation} 
supported in Maslov grading zero. Since the reduced complex is filtered chain homotopy equivalent to $CFK^-(K)$, they have the same Euler characteristic. Equation \eqref{lspacereduced} says in particular that every generator has even Maslov grading, so each contributes a positive term to the Euler characteristic. In other words, if \[ \widetilde{\Delta}_K(T)=\sum_{i=0}^\infty T^{a_i},\] then $\underline{CFK}^-(K)$ has one generator with Alexander grading $a_i$, for each $i\geq 0$. Since multiplication by $U$ is a filtered map (i.e., it never increases the Alexander grading), then necessarily \[M(a_i)=-2i.\]

\begin{figure}
$$\xymatrixcolsep{0.4 pc}\xymatrixrowsep{0.4 pc}
\xymatrix{
& & & & & \underline{A} \\
& & & & \tc{red}{\bu} & 3 & & &&&& & & \bu \ar@{}[]_*{a_0=3}^*{(0)\ }\\
& & & \bu & \bu \ar[l] \ar[dd]& 2\\
& & \bu & \bu \ar[l] \ar[dd] & & 1\\
& \bu & \bu \ar[l] \ar[dd] & & \tc{red}{\bu} & 0 & &&&& & & & \bu \ar@{}[]^*{a_1=0}_*{(-2)\ } \\
\bu & \bu \ar[l] \ar[dd] & & \tc{red}{\bu} & & 1 \ar@{}[]_*{-} & & &&&& & & \bu \ar@{}[]^*{a_2=-1}_*{(-4)\ }\\
& & \bu & & \bu \ar[ll] \ar[d] & 2 \ar@{}[]_*{-}\\
& \bu \dar[dl] & & \bu \ar[ll] \ar[d] & \tc{red}{\bu} & 3 \ar@{}[]_*{-} & &&&& & & & \bu \ar@{}[]^*{a_3=-3}_*{(-6)\ } \\
& & & \tc{red}{\bu} & & 4 \ar@{}[]_*{-} & & & &&&& & \bu \ar@{}[]^*{a_4=-4}_*{(-8)\ } \dar[d]\\
& & & & & \tc{white}{i} & & &&&& & & \\
& & & \ar@{}[]|*{CFK^-(T_{3,4})} & & & & &&&& & & \ar@{}[]|*{\underline{CFK}^-(T_{3,4})}
 }
$$
\caption{To the left is the knot Floer complex $CFK^-$ for the (3,4)-torus knot, which, by Proposition \ref{lspacechar}, is determined by its Alexander polynomial. To the right is the reduced complex, which, for any $L$-space knot, is isomorphic to $\F[U]$, supported in Maslov grading zero. Note that the reduced complex on the right has a generator for each generator colored red on the left -- the bottom-most generator of each staircase summand.}
\label{T34}
\end{figure}
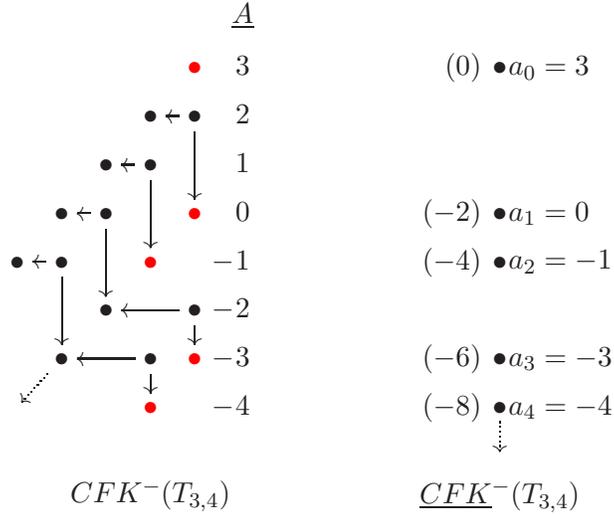

 Figure \ref{T34} gives an illustration for the case of the $(3,4)$-torus knot\footnote{It was shown by Moser \cite{Moser} that torus knots admit lens space (hence $L$-space) surgeries. }, where
\begin{equation}
\Delta_{T_{3,4}}(T) = 1 -(T^2+T^{-2}) + (T^3+T^{-3}),
\label{34poly}
\end{equation}
and therefore 
\begin{equation}
\widetilde{\Delta}_{T_{3,4}}(T) = T^3 +1 +T^{-1}+T^{-3}+T^{-4} +\cdots.
\label{34modpoly}
\end{equation}

\begin{proof}[\sc Proof of Theorem \ref{alternate}]
First note that it is sufficient to prove the propositon for positive surgeries (see Remark \ref{mirror}). 

So, let $K$ be an $L$-space knot in  $S^3$, so that \[\widetilde{\Delta}_K(T)=\sum_{i=0}^\infty T^{a_i}.\] Then the reduced complex $\underline{CFK}^-(K)$ consists of a single generator of Alexander grading $a_i$ and Maslov grading $-2i$, for each $i\geq 0$. Since the $a_i$ are strictly decreasing, it follows that $\widetilde{d}$, as defined in \eqref{shiftddef}, is given by \[\widetilde{d}(K,m) = -2 \min \{ i | a_i\leq m\},\] and therefore 
\begin{equation}
\widetilde{d}(K,a_i-1) = -2 (i+1).
\end{equation} 
Substituting these values into inequality \eqref{shiftdineq}, we obtain
\begin{equation}
i+1\leq \left\lceil \frac{g(K)-(a_i-1)}{2}\right\rceil,
\end{equation}
so
\[ a_i \leq g(K)-2i.\]

\end{proof}
\begin{proof}[\sc Proof of Theorem \ref{newrestriction}]
\begin{figure}
$$\xymatrixcolsep{0.3 pc}\xymatrixrowsep{0.3 pc}
\xymatrix{
& & & & & & & & \underline{A}\\
& & & & & & & \tc{red}{\bu} & 6  & & & && & & & \bu \ar@{}[]_{(0)}\\
& & & & & & \bu & \bu \ar[l] \ar[ddd]^{r_3} & 5\ar@{}[]^*{\ =m_1} & & & & \\
& & & & & \bu & \bu \ar[l] \ar[ddd] & & 4\\
& & & & \bu & \bu \ar[l] \ar[ddd] & & & 3 \\
& & & \bu & \bu \ar[l] \ar[ddd] & & & \tc{red}{\bu} & 2 & && & & & & & \bu \ar@{}[]^{(-2r_1)} \\
& & \bu & \bu \ar[l] \ar[ddd] & & & \tc{red}{\bu} & & 1 & & & & && & & \bu \\
& \bu & \bu \ar[l]_{r_1=1} \ar[ddd] & & & \bu & & \bu \ar[ll] \ar[dd]^{r_2} & 0\ar@{}[]^*{\ =m_2} & & & &   \\
\bu & \bu \ar[l] \ar[ddd] & & & \bu & &  \bu \ar[ll] \ar[dd] & & -1\\
& & & \bu & &  \bu \ar[ll] \ar[dd] & & \tc{red}{\bu} & -2 && & & & & & & \bu \ar@{}[]^{(-2(r_1+r_2))} \\
& & \bu & &   \bu \ar[ll]_{r_2=2} \ar[dd] & & \tc{red}{\bu} & & -3 & & & & && & & \bu \\
& \bu & &  \bu \ar[ll] \ar[dd] & & \tc{red}{\bu} & & & -4 & & & & & & &&  \bu \\
& & \da[ddll] & & \bu & & & \bu \ar[lll]_{r_3=3} \ar[d]^{r_1} & -5\ar@{}[]^*{\ \ =m_3} & & & &  \\
& & & \bu & & & \bu \ar[lll] \ar[d] & \tc{red}{\bu} & -6 & & && & & & & \bu \ar@{}[]^{(-2(r_1+r_2+r_3))} \\
& & & & & & \tc{red}{\bu} & & -7 & & & & & && & \bu \\
& & & &  \ar@{}[]|*{CFK^-(T_{4,5})} & & & & & & & && & & &  \ar@{}[]|*{\underline{CFK}^-(T_{4,5})}
 }
$$
\caption{The complex $CFK^-$ for the $(4,5)$-torus knot, and its reduced form $\underline{CFK}^-$. Note that the integers $r_1, r_2$ and $r_3$ are the horizontal lengths of each staircase, from left to right (and by symmetry, the vertical lengths, from bottom to top). This figure illustrates how the $m_j$'s -- the Alexander gradings at which the reduced complex ``jumps" -- are determined by the $r_i$'s, and further, how the values of $\widetilde{d}(K,m_j)$, given in parentheses to the right, are determined by the $r_i$'s. }
\label{T45}
\end{figure}
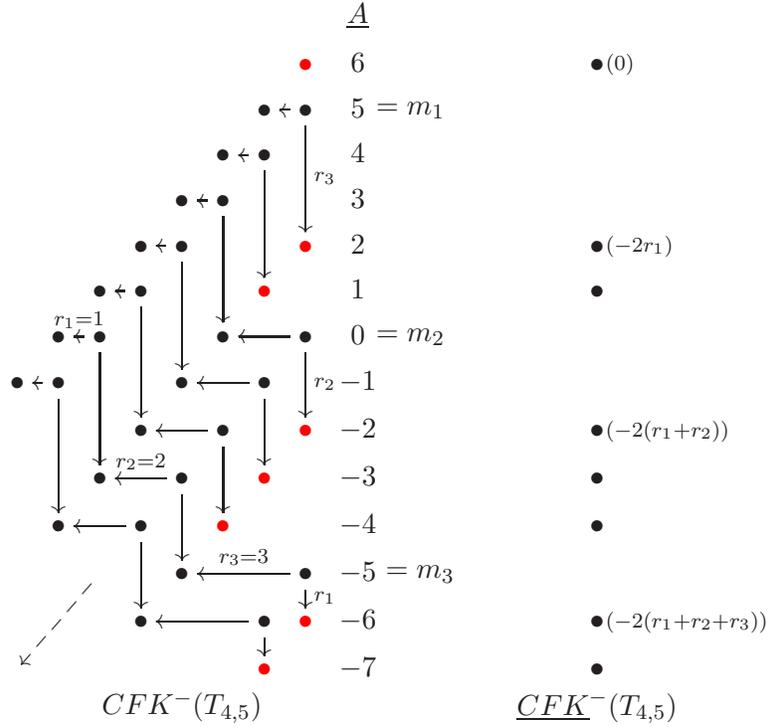

Let $K$ be an $L$-space knot, so that \[\Delta_K(T)=\sum_{i=-k}^{k} (-1)^{k+i}T^{n_i}.\] We have introduced the variables \[ \vec{r}=(r_1, \ldots, r_k ) \] as the ``gaps" in the Alexander polynomial (the difference in the exponents of consecutive nonzero terms), \[r_i=n_{k+2-2i}-n_{k+1-2i}.\] While $\vec{r}$ records only every second gap, by the symmetry of $\Delta(T)$, this determines the polynomial uniquely. Diagramatically,  notice that $\vec{r}$ is simply the list of horizontal lenghts of a staircase summand of $CFK^-$, in order from left to right. See Figure \ref{T45} for the example of the $(4,5)$-torus knot, which has \[\Delta_{T(4,5)}(T)=-1 + (T^2+T^{-2}) - (T^5+T^{-5}) + (T^6+T^{-6}).\]

Next we observe how, given $\vec{r}$, to compute both sides of inequality \eqref{shiftdineq} for any $m$, with Figure \ref{T45} as a guide. We will focus on the values labeled $m_j$ in Figure \ref{T45}; in other words, the values where we have the ``jumps" in the reduced complex. More precisely, if we let \[m_j= g(K)- \left( \sum_{i=1}^j r_i + \sum_{i=k-j+2}^k r_i \right),\] we have that \[\widetilde{d}(K,m_j)=-2\sum_{i=1}^j r_i.\] Substituting these values into inequality \eqref{shiftdineq} when $m=m_j$ gives
\begin{equation}
 \sum_{i=1}^j r_i \leq \left\lceil \frac{g(K)-\left(g(K)- \left( \sum_{i=1}^j r_i + \sum_{i=k-j+2}^k r_i \right)\right)}{2}\right\rceil.
\label{rineq}
\end{equation}
 The case $j=1$ gives \( r_1\leq \lceil \frac{r_1}{2}\rceil,\) so, since each $r_i$ is a positive integer, $r_1$ must equal 1. Substituting this into \eqref{rineq} gives
\begin{align*}
 1+\sum_{i=2}^jr_i \leq & \left\lceil \frac{ 1+ \sum_{i=2}^j r_i + \sum_{i=k-j+2}^k r_i }{2}\right\rceil\\
 \sum_{i=2}^jr_i \leq & \left\lceil \frac{ -1+ \sum_{i=2}^j r_i + \sum_{i=k-j+2}^k r_i }{2}\right\rceil,
\end{align*}
from which it follows that \[ \sum_{i=2}^jr_i \leq \sum_{i=k-j+2}^kr_i .\]
This is sufficient to prove the claim. We could similarly obtain inequalities by considering values of $m$ different from the $m_j$, but those would be no stronger, and therefore provide no more restrictions on $\vec{r}$.
\end{proof}

As an example, consider a knot $K$ with \[ \Delta_{K}(T)= -1 + (T^2+T^{-2}) - (T^3+T^{-3}) + (T^4+T^{-4}), \] so that
\[\widetilde{\Delta}_K(T)=T^4+T^2+T+T^{-2}+T^{-4}+T^{-5}+\cdots. \] This polynomial satisfies the restriction of Proposition \ref{HWrestriction}, but if $K$ were an $L$-space knot, we would have $g(K)=4$, and $a_2=1$. This violates inequality \eqref{aiineq}, so $K$ (and its mirror image) cannot admit an $L$-space surgery. Alternatively, this polynomial has gaps $\vec{r}=\{ 1,2,1\}$, and since $r_2 \nleq r_3$, this violates inequality \eqref{riineq}.

In fact, this completely determines which Alexander polynomials are realized by $L$-space knots of genus less than or equal to 4. All other polynomials satisfying Proposition \ref{HWrestriction} are realized by known $L$-space knots.\\
\begin{center}
\setlength{\tabcolsep}{20pt}
\begin{tabular}{l l}
$\vec{r}$ & $L$-space knot with corresponding $\Delta(T)$ \\ \hline
\{1\} & $T(2,3)$ \\
\{1,1\} & $T(2,5)$\\
\{1,1,1\} & $T(2,7)$\\
\{1,2\} & $T(3,4),$ $(2,3)$-cable of $T(2,3)$\\
\{1,1,1,1\} & $T(2,9)$\\
\{1,1,2\} & $T(3,5)$\\
\{1,2,1\} & excluded by Theorem \ref{newrestriction}\\
\{1,3\} & $(2,5)$-cable of $T(2,3)$
\end{tabular}\\
\end{center}
For knots of genus 5, Theorem \ref{newrestriction} eliminates the polynomials corresponding to $\vec{r}=\{ 1,2,1,1\}$ and $\vec{r}=\{1,3,1\}$, but there are still three more which are not realized by any $L$-space knot known to the author (corresponding to $\{1,1,2,1\}, \{1,2,2\}$ and $\{1,4\}$).


\bibliography{biblens2}{}
\bibliographystyle{amsplain}

\end{document}